\documentclass{aptpub}
\RequirePackage[OT1]{fontenc}
\usepackage[english]{babel}
\usepackage[utf8]{inputenc} 
\usepackage{latexsym,amssymb,amsmath,epsfig,psfrag}

\def\RR{\mathbb{R}}
\def\EE{\mathbb{E}}
\def\NN{\mathbb{N}}

\def\Ent{\mathbf{Ent}}

\def\II{\mbox{ 1\hskip -.29em I}}
\newcommand{\sfrac}[2]{\kern.1em
        \raise.5ex\hbox{$#1$}\kern-.1em
        /\kern-.15em\lower.25ex\hbox{$#2$}}

\newcommand{\eps}{\varepsilon}

\numberwithin{equation}{section}
\authornames{R. Rossignol}
\shorttitle{Arbitrary symmetric threshold widths}

\begin{document}

\title{Arbitrary threshold widths \\for monotone symmetric properties}

\authorone[Universit\'e de Neuch\^atel]{Rapha\"el Rossignol}
\addressone{\\ Facult\'e des Sciences\\ Institut de Math\'ematiques\\ 11 rue
  Emile Argand\\ 2000 Neuch\^atel, SUISSE\\
{\sc e-mail}: raphael.rossignol@unine.ch\\
{\sc URL}: www.math-info.univ-paris5.fr/{\~{\hspace{0mm}}}rost
}

\begin{abstract}  We investigate the threshold widths of some symmetric
 properties which range asymptotically between $1/\sqrt{n}$ and $1/\log n$. These
  properties are built using a combination of failure sets arising from
  reliability theory. This combination of sets is simply called a product. Some general results on the threshold
  width of the product of two sets $A$ and $B$ in terms of the threshold
  locations and widths of $A$ and $B$ are provided.
\end{abstract}

\keywords{threshold width; zero-one law; parallel-series system;
  series-parallel system; $k$-out-of-$n$ system}

\ams{60F20}{60E15, 60K10}

%

\section{Introduction}
Let $n$ be a positive integer, $p$ a real number
in $[0,1]$, and denote by $\mu_{n,p}$ the probability measure on
$\{0,1\}^n$ which is the product of $n$ Bernoulli measures with
parameter $p$.
$$\forall
x\in\{0,1\}^n,\;\mu_{n,p}(x)=p^{\sum_{i=1}^nx_i}(1-p)^{\sum_{i=1}^n(1-x_i)}\;.$$
We write $\mu_p$ instead of $\mu_{n,p}$ when no confusion is
possible. If $A$ is a subset of $\{0,1\}^n$, we say that $A$ is \emph{monotone} if and only if:
$$\left(x\in A\mbox{ and }x\preceq y\right)\Longrightarrow y\in A\;,$$
where $\preceq$ is the partial order on $\{0,1\}^n$ defined
coordinate-wise. It follows from an elementary coupling device that
for $A$ a monotone subset, the mapping  $p\mapsto\mu_p(A)$ is increasing. For many examples of interest (see section
\ref{examples} for some examples), a threshold
phenomenon occurs for property $A$ in the sense that the function
$p\mapsto \mu_p(A)$ ``jumps'' from near 0 to near 1 over a very short interval
of values of $p$. Such threshold phenomena have been shown to occur in most
discrete probabilistic models, such as random graphs (see
Bollob\'as\cite{Bollobas}), percolation (see Grimmett \cite{Grimmett}),
satisfiability in random constraint models (see Creignou and Daud\'e
\cite{CreignouDaude}, Friedgut \cite{Friedgut1}, Bolllob\'as et al. \cite{Bollobasetal01}), local properties in random
images (see Coupier et al. \cite{CoupierDesolneuxYcart05}), reliability (see
Paroissin and Ycart \cite{ParoissinYcart03}) and so on. To make the statement of a
threshold phenomenon more precise, one need first to define the
\emph{threshold width} of a \emph{non trivial} monotone subset $A$. We say that $A$ is
non trivial if it is non empty and different
from $\{0,1\}^n$ itself. When $A$ is non trivial and monotone, the mapping
$p\mapsto\mu_p(A)$ is invertible. Thus, for  $\alpha \in [0,1]$, let  $p(\alpha)$ be the
unique real in $[0,1]$ such that $\mu_{p(\alpha)}(A)=\alpha$. The
threshold width of a subset is the length of the ``transition interval'', that is to say, the interval over which
its probability raises from $\eps$ to $1-\eps$.
\begin{defn}
\label{defwidth}
Let $A$ be a non trivial monotone subset of $\{0,1\}^n$. Let
$\eps\in\left\rbrack 0, 1/2\right\rbrack$. The \emph{threshold width
of $A$ at level $\eps$} is:
$$\tau(A,\eps)=p(1-\eps)-p(\eps)\;.$$
\end{defn}
\begin{figure}[!ht]
\begin{center}
\psfrag{a}[][]{$p(\eps)$}
\psfrag{c}[][]{$p(1-\eps)$}
\psfrag{d}[][]{$\eps$}
\psfrag{e}[][]{$1-\eps$}
\psfrag{f}[][]{$p$}
\psfrag{g}[][]{$\mu_p(A)$}
\psfrag{h}[][]{$\underbrace{\hspace{1cm}}_{\tau(A,\eps)}$}
\includegraphics[width=6cm]{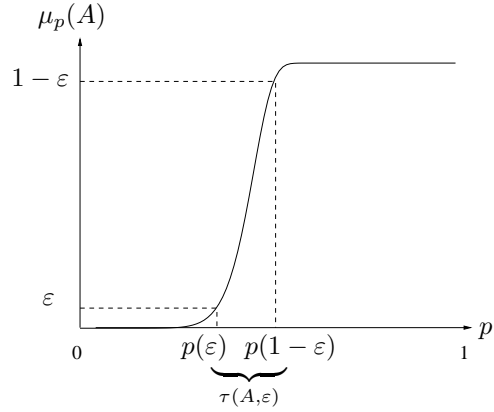}
\end{center}
\caption{Example of a threshold width of level $\eps$.}
\label{fig:seuildef}
\end{figure}
When one investigate the threshold of a monotone property, for example connectivity in the
random graph, one has to do with a sequence of non trivial monotone subsets
$A=(A_n)_{n\in\NN^*}\in\left(\{0,1\}^{\alpha_n}\right)^{\NN^*}$ where
$(\alpha_n)_{n\in\NN^*}$ is an increasing sequence of integers. In the sequel,
we shall suppose that $(\alpha_n)_{n\in\NN^*}$ is only nondecreasing, for
technical reasons. Remark that, in order to get an intrinsec notion of width or localisation
order, one has to keep in mind the size $\alpha_n$ in which the subsets $A_n$
take place. Therefore, if the threshold width of a subset $A_n$ of
$\{0,1\}^{\alpha_n}$ is of order $a(n)$, we should rather express it as
$a\circ\alpha^{-1}(\alpha_n)$, where $\alpha^{-1}$ is the pseudo-inverse of $\alpha$:
$$\forall n\geq\alpha_0,\;\alpha^{-1}(n)=\sup\{k\in\NN\mbox{ s.t. }\alpha_k\leq n\}\;.$$
In order to describe the asymptotic behaviour of a
property, we shall therefore use the following definitions.
\begin{defn}
Let $A=(A_n)_{n\in\NN^*}\in\left(\{0,1\}^{\alpha_n}\right)^{\NN^*}$ be a
monotone property, $a(n)$ and $b(n)$ be two sequences of real numbers in
$[0,1]$, and $\alpha\in [0,1]$.\\
The property $A$ is \emph{ located at
$\alpha$} if:
$$\forall \eps \in ]0,1[,\;p_{A_n,\eps}\xrightarrow[n\rightarrow
+\infty]{}\alpha\;.$$
\emph{ The location of $A$ is of order $a$} if:
$$\forall \eps \in ]0,1[,\;p_{A_n,\eps}=O\left(a(\alpha_n)\right)\;,$$
as $n$ tends to infinity.\\
\emph{ The  threshold width of $A$ is of order $b$} if:
$$\forall \eps \in ]0,1[,\;\tau(A_n,\eps)=O\left(b(\alpha_n)\right)\;,$$
as $n$ tends to infinity.\\
The property $A$ \emph{ has a sharp threshold} if:
$$\forall \eps \in ]0,1[,\;\frac{\tau(A_n,\eps)}{p_{A_n,1/2}(1-p_{A_n,1/2})}\xrightarrow[n\rightarrow
+\infty]{}0\;.$$
The property $A$ \emph{ has a coarse threshold} if it does not have a sharp threshold.
\end{defn}
Intuitively one would be tempted
to say that a subset $A$ will have a narrow threshold unless a few
coordinates have a strong influence on its definition (as an example, think of $A=\{x \mbox{ s.t. }x(1)=1\}$). In many examples,
this idea is captured by the notion of symmetry.
\begin{defn}
\label{defsymmetric}
The subset $A$ of $\{0,1\}^n$ is said to be symmetric if and only if there exists a
subgroup $G$ of $\mathcal{S}_n$ (group of permutations) acting  transitively on
$\{1,\ldots ,n\}$, such that $A$ is invariant under the action of
$G$:
$$\forall g\in G,\;\forall x\in A,\; g.x=\left(x_{g^{-1}(1)},\ldots
  ,x_{g^{-1}(n)}\right)\in A\;.$$
\end{defn}
For a symmetric subset, no coordinate has a stronger influence than any
other. In  Friedgut and Kalai \cite{FriedgutKalai}, it is proven that the threshold width of any
symmetric subset $A\subset\{0,1\}^n$ is at most of order $1/\log n$. For properties whose
threshold is located away from 0 and 1, Friedgut and Kalai show that this upper bound is
tight and that the threshold width
is at least of order $1/\sqrt{n}$. In order to deepen the link between the
invariance group of $A$ and the largest possible threshold width for $A$,
 Bourgain and Kalai \cite{BourgainKalai} introduce, for any permutation group
$G\subset\mathcal{S}_n$,
$$T_G(n,\eps)=\sup\{\tau(A,\eps)\mbox{ s.t. } A \mbox{ is invariant under the
  action of }G\}\;.$$
They obtain nearly optimal asymptotics for $T_G(n,\eps)$ when $G$ is a
  primitive permutation group. Recall that a permutation group
  $G\subset\mathcal{S}_n$ is primitive if its action on $\{1,\ldots ,n\}$ has
  no nontrivial group blocks, where a group block is a subset $B$ of
  $\{1,\ldots ,n\}$ such that for all $g\in G$, $g(B)=B$ or $g(B)\cap
  B=\emptyset$. Essentially, Bourgain and Kalai \cite{BourgainKalai} show that there are some
  gaps in the possible behaviours of $T_G(n,\eps)$ for primitive groups
  $G$. When $G$ is  $\mathcal{A}_n$ (the alternating
  group) or $\mathcal{S}_n$, $T_G(n,\eps)$ is of order $1/\sqrt{n}$. For any
  other primitive group, $T_G(n,\eps)$  is either of order $\log^{-c} n$, for $c$
  belonging to arbitrarily small intervals around a value of the form
  $(k+1)/k$, where $k$ is a positive integer depending only on $G$, or of
  order $\log^{-c(n)} n$, with $c(n)$ which tends to one as $n$ tends to
  infinity. These results concern the worst threshold intervals for a given
  transitive group. In order to complete these results, it is natural to ask, given an increasing sequence of positive real
  numbers $a(n)$ between $\log n$ and $n^{1/2}$, whether there exists a symmetric property
  $A$ whose threshold width is $1/a(n)$. Only few types of such asymptotics
  are known. The main result of this paper is Theorem \ref{thmlargeursvariees}, which gives a
  positive answer to this question under a mild hypothese of smoothness on
  the sequence $a(n)$. This result is achieved by using a
  combination of two properties $A$ and $B$ that we shall simply call the
  product of $A$ and $B$.

This paper is organized as follows. Section \ref{examples} is devoted to some
examples of properties with explicit threshold widths and locations. Some of
them, which arise from reliability theory, will be used further as elementary building blocks to derive more general
widths. In section \ref{product}, we derive the
basic properties of the product of $A$ and $B$ which turns out to have a simple
interpretation in terms of failure sets. We prove that the product of $A$ and $B$ has a threshold
width which is the product of that of $A$ and $B$ as soon as the threshold of
$B$ is located away from 0 and 1. This result allows us to obtain in section
\ref{seuilsintermediaires} some symmetric properties of $\{0,1\}^n$ with arbitrary
threshold widths between $1/\log n$ and $1/\sqrt{n}$. For the sake of
completeness, we also study the case where the threshold of $B$ tends to 0 or
1. Although we do not give an extensive understanding of what may happen, we show in section \ref{coarsesharp} that if $A$ and $B$ have a threshold located respectively in 0 and 1, then $A\otimes B$ has a sharp threshold.

\section{Examples of explicit threshold widths and locations}
\label{examples}
In presenting the following examples of thresholds, our aim is twofold. First,
we want to describe some of the few already known types of behaviour. Second, we shall use
some of
these examples in section \ref{seuilsintermediaires}, to derive more general
widths thanks to the product of properties.

 One of the typical examples of threshold phenomena is that of the random
graphs $\mathcal{G}(n,p(n))$ (see Erd\H{o}s and R\'enyi \cite{ErdosRenyi},
Bollob\'as \cite{Bollobas}, Spencer \cite{Spencer}). The graph $\mathcal{G}(n,p)$ has $n$ vertices,
and each one of the $N=n(n-1)/2$ possible edges is present with probability $p$,
independently from the others. Once a
labelling of the vertices is choosen, one can
denote by $(\{0,1\}^N,\mu_{N,p})$ the probability space of the random graph
$\mathcal{G}(n,p)$.  
\begin{ex} Small balanced subgraphs\\
{\rm Let $H$ be a fixed connected graph with $v$ vertices and $e$ edges, and
  suppose that $H$ is balanced, that is to say none of its subgraphs has
  average degree strictly smaller than $H$. Denote by $A_H$ the property for
  a graph to contain at least one  copy of $H$. The threshold of $A_H$ is
  located at $O(n^{-v/e})$, i.e $O(N^{-v/2e})$, and has width of the same
  order (cf. Spencer \cite{Spencer} for instance). This implies that $A_H$ has a
  coarse threshold.

}
\end{ex}

\begin{ex} Connectivity\\
{\rm It is known (see Bollob\'as \cite{Bollobas84}), that the probability for $\mathcal{G}(n,p(n))$ to
be connected goes from $\eps+o(1)$ to $1-\eps+o(1)$ when $p(n)=\log
  n/n+c/n$, and $c$ goes from $\log
\left(1/\log 1/\eps\right)$ to $\log
\left(1/\log 1/(1-\eps)\right)$. In this example, the threshold
is located around $\log n/n$ i.e $\log N /(2\sqrt{2N})$, and its width is of order
$O\left(1/n\right)$, i.e $O\left(1/\sqrt{N}\right)$. Thus, this threshold is sharp.
}
\end{ex}



Let us turn to examples occuring in reliability theory. In this framework, at
  instant $t$, two characteristic quantities of the system are especially important: the reliability,
  that is the probability that there never occured any breakdown before $t$, and the
  non-availability, which is the probability that the system is down at
  instant $t$ (see for
  instance Barlow and Proschan \cite{BarlowProschan}). Of course, these quantities differ if the
  system is repairable. The analysis of the reliability of large systems, for
  instance its asymptotic behaviour, is generally much more difficult than the
  analysis of the non-availablity. We shall only focus on the latter one, but
  want to stress the fact that when one deals with a large system composed of
  repairable Markovian components, it is natural to expect strong similarities between the asymptotics of
  the two quantities (see for example Paroissin and Ycart \cite{ParoissinYcart04}). When $\mathcal{A}$ denotes a system composed of
$n$  binary components, one can describe the states of these  components as a
state in $\{0,1\}^n$, 1 standing for a failed component, and 0 for a
working component. One can therefore associate to $\mathcal{A}$ its
failure subset, which is the subset $A$ of $\{0,1\}^n$ containing all
the configurations of the $n$ components such that the system
$\mathcal{A}$ fails. If we assume that a component is failed independently from
  the others with probability
  $p$, $\mu_{n,p}$ is the distribution of the state of $\mathcal{A}$ in
  $\{0,1\}^n$, and $\mu_{n,p}(A)$ is the non-availability of
  $\mathcal{A}$. It is very natural to assume that the subset $B_n$ is
  monotone (if the system is down, and a component fails, then the system
  remains down). The question of how quickly $\mu_{n,p}(A)$ ``jumps from 0 to
  1'' is of great importance (see Paroissin and Ycart \cite{ParoissinYcart03} for an application of
  the works of Friedgut and Kalai \cite{FriedgutKalai} and Bourgain et al. \cite{BKKKL} in this context). The main result of this article, Theorem \ref{thmlargeursvariees}, relies on
examples \ref{exksurn} and \ref{exparser}.
\begin{ex} \label{exksurn}$k$-out-of-$n$ system\\
{\rm The system is failed when the total number of failed components is
  greater than a certain threshold $k(n)$. The failure subset is therefore:
$$A_{k,n}=\left\lbrace x\in \{0,1\}^n\mbox{ s.t. }\sum_{i=1}^nx_i\geq k\right\rbrace\;.$$
Note that the particular cases of
 $A_{m-1,m}$ and  $A_{1,r}$ correspond respectively to parallel and series
 system. Obviously, $A_{k,n}$ is monotone and invariant under every permutation of the coordinates. It is therefore a monotone symmetric subset of $\{0,1\}^n$. Since the sum $\sum_{i=1}^nx_i$ has mean $np$ and variance $np(1-p)$
when $x$ is distributed according to $\mu_p$, one can guess,
intuitively, that $A_{k,n}$ has a threshold located at $k/n$,
and of order $\sqrt{(k/n)\times (1-k/n)}/\sqrt{n}$. We shall precise this
intuition when $k=\lfloor n/2\rfloor$ in Lemma \ref{lemseuilsqrtn}.
}
\end{ex}

\begin{ex} \label{exparser}Parallel-series system\\
{\rm A parallel-series system contains $n=r\times m$
components which are assembled into
$r$ blocks containing $m$ composants. The system is failed as soon as a block is
failed, and a block fails if all its components are failed. Of course, the
non-availability of such a system is very easy to derive. Let $B_n$ denote its
failure subset:
$$\mu_p(B_n)=1-(1-p^m)^r\;.$$
For example, when $m=\lfloor\log_2 k\rfloor$, $r=\lfloor
k/\log_2k\rfloor$ and $k\geq 2$, the threshold of $B_n$ is
located at $1/2$ with a width of order $1/\log n$ (see Lemma
\ref{lemseuillogn} below). Remark that $B_n$ is monotone and symmetric
(under permutation of the components inside a block and permutation of the
blocks). Such systems, with multi-states components instead of binary ones,
have been studied by
Kolowrocki \cite{Kolowrocki94a,Kolowrocki94b}, and a concrete application is presented in
\cite{Kolowrocki01}. One can also define the dual system called
series-parallel system, in which components are assembled into
$r$ blocks containing $m$ composants, the system is failed when all blocks are
failed, and a block is failed as soon as one component is failed.
}
\end{ex}

\begin{ex} \label{exkconn}Consecutive $k$-out-of-$n$ system\\
{\rm Components are arranged around a circle. The system is failed as soon as
  there are at least $k(n)$ consecutive components down. This model has
  an asymptotic behaviour similar to the Parallel-series system with $\lfloor
  n/k\rfloor$ blocks of $k$ components. For example, when $k=\lfloor n/\log_2
  n\rfloor$, the threshold of the failure subset is located at $1/2$, with a
  width of order $1/\log n$ (for a similar result, see Paroissin and Ycart 
  \cite{ParoissinYcart03}). This model was introduced by Kontoleon \cite{Kontoleon01} to model
  some problems arising in engineering science, such as oil transportation
  using pipelines, telecommunication system by spacecraft relay station or
  transmission of data in a ring of computer ring networks, etc.
}
\end{ex}

\section{The product of subsets of $\{0,1\}^n$}
\label{product}
\begin{figure}
\begin{center}
\psfrag{a}[][]{$\mathcal{A}$:}
\psfrag{b}[][]{$\mathcal{B}$:}
\psfrag{A}[][]{\begin{minipage}[b]{5cm} Set of failure states $A$,\\
    Parallel system.\end{minipage}}
\psfrag{B}[][]{\begin{minipage}[b]{5cm} Set of failure states $B$,\\
    Series system.\end{minipage}}
\psfrag{C}[][]{\begin{minipage}[b]{5cm} Set of failure states
    $A\otimes B$,\\ Parallel-series system.\end{minipage}}
\psfrag{D}[][]{\begin{minipage}[b]{5cm} Set of failure states
    $B\otimes A$,\\ Series-parallel system.\end{minipage}}
\includegraphics[width=10cm]{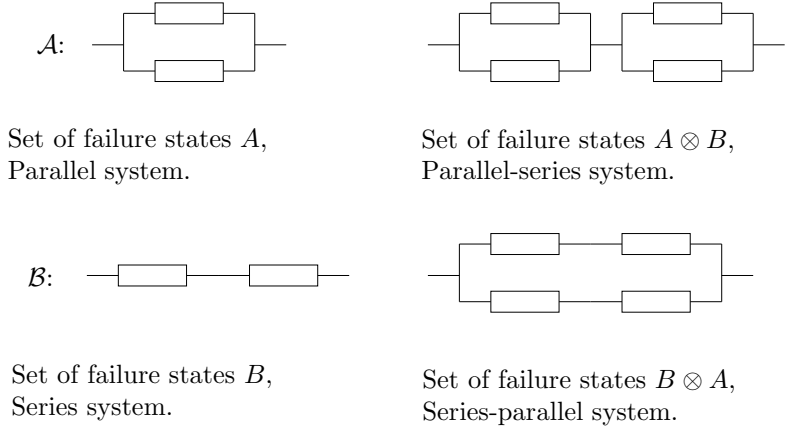}
\caption{Parallel-series and series-parallel systems are obtained via a product.}
\label{fig:tenseurparser}
\end{center}
\end{figure}

\begin{figure}[!ht]
\begin{center}
\psfrag{a}[][]{$\mathcal{A}$:}
\psfrag{b}[][]{$\mathcal{B}$:}
\psfrag{A}[][]{Set of failure states $A$.}
\psfrag{B}[][]{Set of failure states $B$.}
\psfrag{C}[][]{Set of failure states $A\otimes B$.}
\includegraphics[width=11cm]{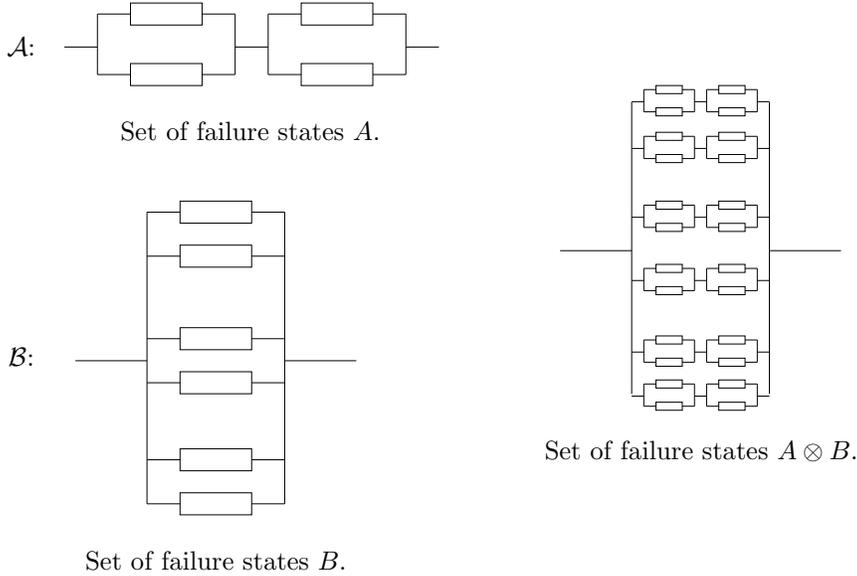}
\caption{An example of product in reliability theory.}
\label{fig:tenseurfiab}
\end{center}
\end{figure}

As far as we know, whereas the influence of simple operations between
properties has been extensively studied whithin the so-called 0-1 laws which
occur in logic (see Compton \cite{Compton}), no such work has been undertaken
regarding the threshold phenomena. The first combinations of properties that
come to mind, union and intersection, behave quite in an unpleasant way with
respect to the threshold width (see \cite{Rossthese}, chapter 3). In this
section, we will show the nice behaviour of another 
combination which we simply call the product. Even though linearity does not
play any role in this setting, it is worth noting the similarity between this product
and the Kronecker product of matrices. Given two properties $A$ and $B$, on two distinct spaces, their product is a property combining the belongings to $A$ and
$B$ in the following way.
\begin{defn}
\label{deftens}
Let $A$  be a subset of $\{0,1\}^r$ and $B$ a subset of
$\subset\{0,1\}^m$. The \emph{ product of $A$
and $B$}, denoted by $A\otimes B$ is the subset of $\left(\{0,1\}^r\right)^m$ defined by:
$$\eta\in A\otimes B\Leftrightarrow \left(\II_{\eta_1\in A},\ldots,\II_{\eta_m\in A}\right)\in B\;,$$
where
$$\eta=(\eta_1,\ldots ,\eta_m)\mbox{ and }\;\forall j\in\{1,\ldots ,m\},\;\eta_j\in\{0,1\}^r\;.$$
\end{defn}
In order to visualize the precise meaning of this definition, it
is convenient to consider this product via the language of
reliability theory. Let $A$ denote the failure set of a system
$\mathcal{A}$ composed of
$r$ components, and  $B$ be the failure set of another
system $\mathcal{B}$, with $m$ components. Then $A\otimes B$ is the failure subset
of the system obtained by replacing the components in $\mathcal{B}$
by $m$ independent copies of $\mathcal{A}$. For example, one can obtain the
so-called  parallel-series and series-parallel systems from some
elementary building blocs: the series and parallel systems (see figure
\ref{fig:tenseurparser}). This building set can be continued,
embedding systems one in another (see figure \ref{fig:tenseurfiab}).

Now, let us describe the basic properties of this product. A very nice feature is the link between the probability
of $A\otimes B$ and those of
$A$ and $B$. It is also easy to get some invariance and monotonicity properties for $A\otimes B$
providing some similar hypotheses for $A$ and $B$. In the sequel, if $\eta=(\eta_1,\ldots ,\eta_m)$ belongs to $\left(\{0,1\}^r\right)^m$, with $\eta_j\in\{O,1\}^r$ for every $j$, we will denote by $\eta_{i,j}$ the $i$-th coordinate of $\eta_j$, which is therefore 0 or 1. In this way, we identify $\left(\{0,1\}^r\right)^m$ and $\{0,1\}^{\{1,\ldots, r\}\times\{1,\ldots,m\}}$.
\begin{prop}
\label{propprodtensgen}
Let $A\subset\{0,1\}^r$ and $B\subset\{0,1\}^m$. 
\begin{enumerate}
\item For every $p$ in $[0,1]$,
$$\mu_{mr,p}(A\otimes B)=\mu_{m,\mu_{r,p}(A)}(B)\;,$$
\item If $A$ and $B$ are monotone, then $A\otimes B$ is monotone.
\item If $A$ is invariant under the action of a subgroup $G$ of
  $\mathcal{S}_r$ and $B$ is invariant under the action of a subgroup
  $H$ of $\mathcal{S}_m$, then $A\otimes B$ is invariant under the
  action of the subgroup $G\times H$ of the permutations of
  $\{1,\ldots, r\}\times\{1,\ldots,m\}$ defined by:
$$\forall i\in\{1,r\},\;\forall j\in\{1,m\},\;(g,h).(i,j)=(g.i,h.j)\;.$$
\end{enumerate}
\end{prop}

\begin{proof}
If $\left(\eta_1,\ldots,\eta_m\right)$ are independant and distributed
according to the law $\mu_{r,p}$, then $\left(\II_{\eta_1\in
    A},\ldots,\II_{\eta_m\in A}\right)$ has law
$\mu_{m,\mu_p(A)}$. This proves the first assertion \\
Let us prove now the second assertion. Let $\eta$ and $\zeta$ belong
to $\left(\{0,1\}^r\right)^m$. Suppose that $\eta\leq\zeta$, i.e $$\forall i\in\{1,\ldots,m\},\;\eta_i\leq\zeta_i\;.$$
Since $A$ is monotone,
\begin{equation}
\label{eqcroisstens1}
\left(\II_{\eta_1\in A},\ldots,\II_{\eta_m\in A}\right)\preceq \left(\II_{\zeta_1\in A},\ldots,\II_{\zeta_m\in A}\right)\;.
\end{equation}
Suppose now that $\eta\in A\otimes B$.
\begin{equation}
\label{eqcroisstens2}
\left(\II_{\eta_1\in A},\ldots,\II_{\eta_m\in A}\right)\in B\;.
\end{equation}
Since $B$ is monotone, il follows from (\ref{eqcroisstens1}) and
(\ref{eqcroisstens2}) that $\zeta\in B$, which proves the monotonicity
of $A\otimes B$.\\
Let us prove now the last point of proposition \ref{propprodtensgen}. Let $\eta\in
A\otimes B$, $(g,h)\in G\times H$ and let us denote
$\zeta=(g,h).\eta$. 
$$\zeta_{i,j}=\eta_{(g,h).(i,j)}=\eta_{g.i,h.j}\;,$$
which can be restated as:
$$\zeta=\left(g.\eta_{h.1},\ldots,g.\eta_{h.m}\right)\;.$$
On the other hand,
$$\eta=\left(\eta_1,\ldots,\eta_m\right)\;,$$ 
with $\eta_i\in\{0,1\}^r$. And also:
$$\left(\II_{\eta_1\in A},\ldots,\II_{\eta_m\in A}\right)\in B\;.$$
Therefore,
$$\left(\II_{g(\eta_1)\in A},\ldots,\II_{g(\eta_m)\in A}\right)\in B\;,$$
$$h.\left(\II_{g(\eta_1)\in A},\ldots,\II_{g(\eta_m)\in A}\right)\in B\;,$$
which means:
$$\left(\II_{g(\eta_{h.1})\in A},\ldots,\II_{g(\eta_{h.m})\in A}\right)\in B\;.$$
Thus $\zeta\in A\otimes B$, and the proof is complete.
{\hfill $\square$ \noindent} 
\end{proof}

Intuitively, the first assertion in Proposition \ref{propprodtensgen} suggests that if the threshold of $B$ is located away from 0 and
1, the threshold effects of $A$ and $B$ will conjugate and give birth to a
threshold width the order of which will be the product of the widths of $A$
and $B$. This is indeed the case, and this is roughly the statement of
Proposition \ref{propprodtensseuil}. Actually, this result is valid as long as the threshold of $B$
is located away from zero and one, and some additional hypotheses of
homogeneity hold for the threshold widths of $A$ and $B$. When a threshold phenomenon occurs for a property $A$, it is usually true that the
threshold width is homogeneous, in the sense that all the transition
intervals shrink at the same speed. This allows to consider the exact order
of the threshold width, since this one does not depend on the
level $\eps$. We will use
the following definitions of homogeneity and strong homogeneity.

\begin{defn}
\label{defseuilhomogene} Let $A\subset\{0,1\}^{\alpha_n}$ be a
non trivial monotone property, and $\left(a_n\right)_{n\in\mathbb{N}}$ be a sequence of positive real
numbers.\\
The threshold width of the property $A$ is  \emph{ homogeneous of order
  $a_n$} if:
$$\forall \beta,\;\gamma \in ]0,1[,\mbox{
  s.t. }\beta<\gamma,\;p_{A,\gamma}-p_{A,\beta}=\Theta(a_n)\;.$$
The threshold width of the property $A$ is  \emph{ strongly homogeneous of order
  $a_n$} if in addition, for all sequences of real numbers
$\left(\beta_n\right)_{n\in\mathbb{N}}$ and
$\left(\gamma_n\right)_{n\in\mathbb{N}}$ such that 
$$\exists\eps\in]0,1[,\;\forall n\in\mathbb{N},\;\eps<\beta_n<\gamma_n< 1-\eps\;,$$
we have
$$p_{A,\gamma_n}-p_{A,\beta_n}=O((\gamma_n-\beta_n)b_n)\;.$$
\end{defn}
We are now able to state the main result about the width of a product.
\begin{prop}
\label{propprodtensseuil}Let $(r_n)_{n\in\NN}$ and $(m_n)_{n\in\NN}$
be two nondecreasing sequences of integers, $A\subset\{0,1\}^{r_n}$
and $B\subset\{0,1\}^{m_n}$ be two monotone properties.
Suppose that the threshold width of $A$ is strongly
homogeneous of
order $a_n$, and the threshold width of $B$ is homogeneous of order
$b_n$. Suppose in addition that the threshold of $B$ is located away from 0
and 1:
$$\forall \eps\in ]0,1[,\;\exists\delta\in]0,1[,\;\forall n\in \NN^*,\;\delta<p_{B,\eps}<1-\delta\;.$$
Then, the threshold of $A\otimes B\subset\{0,1\}^{r_nm_n}$
has a homogeneous width of order $a_n\times b_n$. Moreover, if the threshold of $A$
is located at $\alpha\in [0,1]$, so does the threshold of $A\otimes B$.
\end{prop}
\begin{proof}
Let $\eps$ be a real number in $]0,1/2[$. According to proposition \ref{propprodtensgen}, 
$$\mu_{m_nr_n,p}(A\otimes B)=\mu_{m_n,\mu_{r_n,p}(A)}(B)\;,$$
Therefore,
$$\mu_{r_n,p_{A\otimes B,\eps}}(A)=p_{B,\eps}\;.$$
Then,
$$p_{A\otimes B,\eps}=p_{A,p_{B,\eps}}\;.$$
$$p_{A\otimes B,1-\eps}-p_{A\otimes B,\eps}=p_{A,p_{B,1-\eps}}-p_{A,p_{B,\eps}}\;.$$
Since the threshold width of $B$ is of order $b_n$,
$$p_{B,1-\eps}-p_{B,\eps}=\Theta(b_n)\;,$$
Recall that, by hypothese,
$$\exists\delta\in]0,1[,\;\forall n\in \NN,\;\delta<p_{B,\eps}<p_{B,1-\eps}<1-\delta\;.$$
Thus, the fact that $A$ has a strongly homogeneous threshold width
of order $a_n$ (see definition \ref{defseuilhomogene}) implies that:
$$p_{A,p_{B,1-\eps}}-p_{A,p_{B,\eps}}=\Theta\left(\left(p_{B,1-\eps}-p_{B,\eps}\right)a_n\right)=\Theta(a_nb_n)\;.$$
Therefore, the threshold of $A\otimes B\subset\{0,1\}^{r_nm_n}$
has a homogeneous width of order $a_n\times b_n$.\\
Now, suppose that $A$ is located at $\alpha\in [0,1]$. Let $\eps$ be a real
number in $]0,1[$. Recall that
$$\exists\delta\in]0,1[,\;\forall n\in
\NN,\;\delta<p_{B,\eps}<1-\delta\;.$$
Since $p_{A\otimes B,\eps}=p_{A,p_{B,\eps}}$,
$$\forall n\in
\NN,\;p_{A,\delta}<p_{A\otimes B,\eps}<p_{A,1-\delta}\;.$$
Thus $p_{A\otimes B,\eps}$ tends to $\alpha$ as $n$ tends to infinity. This
completes the proof.
{\hfill $\square$ \noindent} \end{proof}

\section{Symmetric threshold widths between $1/\log n$ and $1/\sqrt{n}$}
\label{seuilsintermediaires}
In this section, we show how to derive from Proposition
\ref{propprodtensseuil} a large variety of threshold widths, ranging
from $1/\log n$ to $1/\sqrt{n}$. To this end, we need some elementary
building blocks the threshold of which are easy to study, and which we shall
eventually combine in order to obtain the desired threshold width. These
blocks will be taken from the reliability examples of section \ref{examples}.\\
Recall that for any $k\in\{1,\ldots,n\}$, we denote by $A_{k,n}$ the following
subset
of configurations in $\{0,1\}^n$ (see example \ref{exksurn}):
$$A_{k,n}=\{x\in \{0,1\}^n\mbox{ s.t. }\sum_{i=1}^nx_i\geq k\}\;.$$
In the sequel, we shall use $A_{\lfloor n/2\rfloor, n}$ and $A_{1,r}\otimes A_{m-1,m}$ for different values of $n$, $r$ and $m$.
\begin{lem}
\label{lemseuilsqrtn}
Let $A=\left(A_{\lfloor   n/2\rfloor,n}\right)_{n\in\NN^*}$. For every $n\in \NN^*$,
$$\tau(A_{\lfloor   n/2\rfloor,n},\eps)\leq 2\sqrt{\log(1/\eps)/(2n)}\;.$$
Moreover, $A$ has a strongly homogeneous threshold, located at $1/2$, with a width of order $1/\sqrt{n}$.
\end{lem}
\begin{proof}
The simplest way to show that $A$ has a trheshold located at $1/2$ with a
width of order $1/\sqrt{n}$ is perhaps to use the concentration property of
the binomial law. Indeed, Hoeffding's inequality \cite{Hoeffding}
ensures that:
\begin{equation}
\label{inegGDbingauche}\forall \lambda >0,\;\mu_p\left(\sum_{k=1}^{n}
x_i-np> \lambda\sqrt n \right) \leq e^{-2\lambda^2} \; ,
\end{equation}
and 
\begin{equation}
\label{inegGDbindroite}\forall \lambda >0,\;\mu_p\left(\sum_{k=1}^{n}
x_i-np< -\lambda\sqrt n \right) \leq e^{-2\lambda^2} \; .
\end{equation}
Let $\eps$ belong to $]0,1[$, and let
$c=\sqrt{\frac{\log(1/\eps)}{2}}$, so that $\exp(-2c^2)=\eps$. 
If  $p(\eps)$ is such that $\mu_{n,p(\eps)}(A_{\lfloor   n/2\rfloor,n})=\eps$,
then $\lfloor n/2 \rfloor $ cannot
be too far away from $np(\eps)$. Inequality (\ref{inegGDbingauche}) and
(\ref{inegGDbindroite}) imply
that:
$$\lfloor n/2\rfloor- \sqrt{\frac{n\log
    \frac{1}{1-\eps}}{2}}\leq p(\eps) \leq \lfloor n/2\rfloor + \sqrt{\frac{n\log
    \frac{1}{1-\eps}}{2}}\;.$$ 
Therefore, the threshold of $A_{\lfloor n/2\rfloor,n}$ is located at $1/2$:
$$\forall \eps\in ]0,1[,\;p(\eps)\xrightarrow[]{}\frac{1}{2}\;,$$
and its threshold width is at most of order $1/\sqrt{n}$:
$$\forall n\in\NN^*,\;\forall \eps\in ]0,1/2[,\;\tau(A_{\lfloor
  n/2\rfloor,n},\eps)\leq 2\sqrt{\frac{\log
    \frac{1}{\eps}}{2n}}\;.$$
To see that this is the right order, one can express the derivative $d\mu_p(A)/dp$ as follows:
\begin{eqnarray*}
\frac{d\mu_p(A)}{dp}&=&\sum_{x\in\{0,1\}^n}\II_A(x)\frac{d\mu_p(x)}{dp}\;,\\
&=&\sum_{i=1}^n\sum_{x\in\{0,1\}^n}\II_A(x)\frac{x_i-p}{p(1-p)}\mu_p(x)\;,\\
&=&\frac{1}{p(1-p)}Cov\left(\II_A,S_n\right)\;.
\end{eqnarray*}
Then, by Cauchy-Schwarz inequality,
\begin{eqnarray}
\nonumber \frac{d\mu_p(A)}{dp}&\leq&\frac{1}{p(1-p)}\sqrt{\mu_p(A)(1-\mu_p(A))}\sqrt{np(1-p)}\;,\\
\nonumber \frac{d\mu_p(A)}{dp}&\leq&\frac{\sqrt{\mu_p(A)(1-\mu_p(A))}}{\sqrt{p(1-p)}}\sqrt{n}\;,\\
\label{ineqdiff} \frac{1}{\sqrt{\mu_p(A)(1-\mu_p(A))}}\frac{d\mu_p(A)}{dp}&\leq&\frac{\sqrt{n}}{\sqrt{p(1-p)}}\;.
\end{eqnarray}
One can easily integrate this differential inequation. Let us define:
$$J(x)=(1-x)\sqrt{x(1-x)}+{\rm Arctan}\left(\sqrt{\frac{x}{1-x}}\right)\;.$$
We have:
$$\forall x\in ]0,1[,\;J'(x)=\frac{1}{\sqrt{x(1-x)}}\;.$$
Therefore, integrating \ref{ineqdiff} between $p(\eps)$ and $p(1-\eps)$ gives:
$$J(1-\eps)-J(\eps)\leq\sqrt{n}\left(J(p(1-\eps))-J(p(\eps))\right)\;.$$
When $n$ tends to infinity, $p(\eps)$ and $p(1-\eps)$ tend to $1/2$. Therefore:
\begin{eqnarray*}
J(1-\eps)-J(\eps)&\leq&\sqrt{n}J'(1/2)\left((p(1-\eps)-p(\eps))+o(p(1-\eps)-p(\eps))\right)\;,\\
p(1-\eps)-p(\eps)&\leq &\frac{J(1-\eps)-J(\eps)}{2\sqrt{n}}+o\left(\frac{1}{\sqrt{n}}\right)\;.
\end{eqnarray*}
Therefore, $\tau(A,\eps)=O(1/\sqrt{n})$.

Finally, to prove the strong homogeneity of the width, we need a sharp
minoration of $d\mu_p(A)/dp$. A smooth way to do this is to use one of the
discrete isoperimetric inequalities à la Margulis-Talagrand. The work of
Margulis \cite{Margulis} has impulsed a number of more and more accurate
discrete isoperimetric inequalities (see 
\cite{Talagrand93,Bobkov97,BobkovGoetze99,TillichZemor00}). For example let
$\phi$ denote the Gaussian density, i.e
$\phi(t)=\frac{1}{\sqrt{2\pi}}e^{-t^2/2}$, and $\Phi$ the Gaussian cumulative
distribution, i.e $\Phi(x)=\int_{-\infty}^x\phi(t)\;dt$. One can deduce from
the main result in Tillich and Z\'emor \cite{TillichZemor00} that: 
\begin{equation}
\label{ZT}
\frac{d\mu_p(A)}{dp}\geq \frac{\sqrt{n}}{p\sqrt{\log 1/p}}\Psi(\mu_p(A))\;,
\end{equation}
where $\Psi$ stands for $\phi \circ \Phi^{-1}$. Therefore,
$$\frac{d\mu_p(A)}{dp}\geq \sqrt{ne}\Psi(\mu_p(A))\;,$$
Let $\left(\beta_n\right)_{n\in\mathbb{N}}$ and
$\left(\gamma_n\right)_{n\in\mathbb{N}}$ be two sequences of real numbers in $]0,1[$. Integrating this inequality between $p(\beta_n)$ and $p(\gamma_n)$ leads to:
$$p(\gamma_n)-p(\beta_n)\leq\frac{1}{\sqrt{ne}}\int_{\beta_n}^{\gamma_n}\Psi(u)\;du\;$$
Now, suppose that  there exists $\eps \in]0,1[$ such that 
$$\forall n\in\mathbb{N},\;\eps<\beta_n<\gamma_n< 1-\eps\;.$$
Since $\Psi$ is continuous and strictly positive on $]0,1[$,
$$\int_{\beta_n}^{\gamma_n}\Psi(u)\;du=\Theta(\gamma_n-\beta_n)\;.$$
Finally,
$$p(\gamma_n)-p(\beta_n)=O((\gamma_n-\beta_n)/\sqrt{n})\;,$$
and thus, the strong homogeneity of the threshold width of $A$ holds.
{\hfill $\square$ \noindent} \end{proof}

Now, consider a Parallel-series system composed of $r$ blocks containing $m$
components (see example \ref{exparser}). The system is failed as soon as a block is
failed, and a block fails if all its components are failed. The failure subset of such a system is $A_{m,m}\otimes A_{1,r}$. It is symmetric and
monotone. One can easily derive explicitly its probability:
$$\mu_p(A_{m,m}\otimes A_{1,r})=1-(1-p^m)^{r}\;.$$
Thus, for any $\alpha\in ]0,1[$,
$$p_\alpha=\left(1-(1-\alpha)^{1/r}\right)^{1/m}\;.$$
\begin{lem}
\label{lemseuillogn}
For every $k\in\NN^*$, let $B_k$ be the following Parallel-series failure subset:
$$B_k=A_{\lfloor\log_2 k\rfloor,\lfloor\log_2 k\rfloor}\otimes A_{1,\lfloor k/\log_2k\rfloor}\subset\{0,1\}^K\;,$$
where $\log_2$ stands for the logarithm to base 2, and $K=\lfloor\log_2 k\rfloor\times\lfloor k/\log_2k\rfloor$. The property $B_k$ has a sharp threshold located at $1/2$, with a homogeneous width of order $1/\log K$. More precisely, the threshold width of $B_k$ satisfies the following asymptotic expansion:
$$\tau(B_k,\eps)=\frac{\log\frac{\log\frac{1}{\eps}}{\log\frac{1}{1-\eps}}}{2\log K}+o\left(\frac{1}{\log_2 K}\right)\;.$$
\end{lem}
\begin{proof}
Let $\alpha$ belong to $]0,1[$. We have:
\begin{eqnarray*}
p_\alpha&=&\left(1-(1-\alpha)^{1/\lfloor k/\log_2k\rfloor}\right)^{1/\lfloor\log_2 k\rfloor}\;,\\
&=&\left(\frac{\log_2k}{k}\log\frac{1}{1-\eps}+o\left(\frac{\log_2k}{k}\right)\right)^{\frac{1}{\log_2 k}+o\left(\frac{1}{\log k}\right)}\;,\\
&=&\exp\left(-\frac{\log k}{\log_2k}+\frac{\log\log_2 k}{\log_2k}+\frac{\log\log \frac{1}{1-\eps}}{\log_2 k}+o\left(\frac{1}{\log k}\right)\right)\;,\\
&=&\frac{1}{2}\left(1+\frac{\log\log_2 k}{\log_2k}+\frac{\log\log \frac{1}{1-\eps}}{\log_2 k}+o\left(\frac{1}{\log k}\right)\right)\;.
\end{eqnarray*}
Therefore, the threshold of $B_k$ is located at $1/2$. Moreover, for any $\eps\in ]0,1/2[$,
$$\tau(B_k,\eps)=\frac{\log 2\log\frac{\log\frac{1}{\eps}}{\log\frac{1}{1-\eps}}}{2\log k}+o\left(\frac{1}{\log k}\right)\;.$$
Since $K(k)=\lfloor\log_2 k\rfloor\times\lfloor k/\log_2k\rfloor$ is equivalent to $k$ as $k$  tends to infinity, 
$$\tau(B_k,\eps)=\frac{\log 2\log\frac{\log\frac{1}{\eps}}{\log\frac{1}{1-\eps}}}{2\log K}+o\left(\frac{1}{\log K}\right)\;.$$
Thus, $B_k$ has a homogeneous width of order $1/\log K$. 
{\hfill $\square$ \noindent} \end{proof}
We are now able to prove our main result.
\begin{thm}
\label{thmlargeursvariees}
Let $c(n)$ be a nondecreasing sequence of integers such that:
$$\forall n\in\NN^*,\;\log n\leq c(n) \leq \sqrt{n}\;,$$
and suppose that:
$$c(2n)=\Theta\left(c(n)\right)\;.$$
Then, there exists an increasing sequence of integers $(N(n))_{n\in \NN^*}$ and a monotone symmetric property $C\subset\{0,1\}^{N(n)}$ whose threshold
is located at $1/2$, with a width of order $1/c(N(n))$.
\end{thm}
\begin{rem}
{\rm The condition $c(2n)=\Theta\left(c(n)\right)$ implies that $c$ increases
   rather smoothly. Because of the way we choosed to build the property
   mentioned in Theorem \ref{thmlargeursvariees}, we
   cannot avoid this condition. Of course, any condition of the type
   $c(rn)=\Theta\left(c(n)\right)$, with $r$ an integer strictly greater than 1
   would be sufficient, since $c$ is nondecreasing. Remark that in most
   natural choices of a nondecreasing function, this
   condition is fulfilled. Nevertheless, one can build some ``unnatural''
   examples where this condition fails. Indeed, let $(a_j)_{j\in \NN^*}$ be the
   increasing sequence of integers defined
   recursively by:
$$\left\lbrace \begin{array}{l}a_0=2\\ a_{j+1}=e^{\lfloor
      \sqrt{2a_j}\rfloor}\end{array}\right.$$
Remark that $a_{j+1}-1\geq 2a_j$ for every $j\geq 1$. Let $c(n)$ be the nondecreasing sequence of integers defined by:
\begin{eqnarray*}
\forall n\in [a_j,2a_j-1],\;c(n)=\log(n)\;,\\
\forall n\in [2a_j,a_{j+1}-1],\;c(n)=\sqrt{2a_j}\;.
\end{eqnarray*}
On one hand,
$$\forall n\in \NN^*,\;\log(n)\leq c(n)\leq \sqrt{n}\;,$$
and on the other hand,
$$\forall j\in\NN^*,\;c(a_j)=\log a_j\mbox{ and }c(2a_j)=\sqrt{2a_j}\;.$$
Therefore, $c(2n)/c(n)$ is not upperbounded.
}
\end{rem}
\vskip 2mm\noindent {\it Proof of Theorem \ref{thmlargeursvariees}} : 
For any integer $k\geq 2$, let $B_k$ be the failure subset of a
parallel-series system composed of $\lfloor k/\log_2 k\rfloor$
blocks ocontaining $\lfloor \log_2 k\rfloor$ components. Suppose that  $1\leq a(n)\leq n$. Let us define the following integer:
$$N=\left\lfloor \frac{n}{a(n)}\right\rfloor\times a(n)\;,$$
and the following monotone
symmetric subset of $\{0,1\}^N$:
$$C_N=A_{\lfloor a(n)/2\rfloor,a(n)}\otimes B_{\left\lfloor
    n/a(n)\right\rfloor}\;.$$
According to Lemma \ref{lemseuilsqrtn}, $A_{\lfloor a(n)/2\rfloor,a(n)}$ has a threshold located at $1/2$ with a strongly homogeneous width of order $1/\sqrt{a(n)}$. From Lemma \ref{lemseuillogn}, $B_{\left\lfloor
    n/a(n)\right\rfloor}$ has a threshold located at $1/2$ with a homogeneous width of order
$1/\log \lfloor   n/a(n)\rfloor$. Therefore, one can deduce from Proposition \ref{propprodtensseuil}
that $C_N$ has a threshold located at $1/2$, with a homogeneous width of order
$\sqrt{a(n)}\times\log \lfloor   n/a(n)\rfloor$.  Now, one can see that
this device allows us to get nearly any order of threshold width between $1/\log N$ and
$1/\sqrt{N}$. Indeed, let  $\phi_n$ be the following function:
$$\phi_n:\left\lbrace \begin{array}{lcl}[1,n]&\rightarrow &\RR \\
    x&\mapsto&x\left(\log\frac{n}{x}\right)^2\end{array}\right.$$
The derivative of $\phi_n$ is easy to compute:
$$\forall x\in [1,n],\;\phi_n'(x)=\log\frac{n}{x}\times\log
\frac{ne^{-2}}{x}\;.$$
Therefore, $\phi_n$ is a bijection from $[1,n/e^2]$ to $[(\log
n)^2,4n/e^2]$. Let $c(n)$ be a sequence of integers such that:
$$\forall n\in\NN^*,\;\log n\leq c(n) \leq \sqrt{n}\;.$$
Define $\tilde{c}(n)=\inf\lbrace c(n),2\sqrt{n}/e\rbrace$. Thus,
$$\forall n\in\NN^*,\;(\log n)^2\leq (\tilde{c}(n))^2 \leq 4n/e^2\;.$$
Let $a(n)=\phi_n^{-1}\left((\tilde{c}(n))^2\right)\in [1,n/e^2]$. The subset $C_N$ has a threshold width
of order $1/\sqrt{\phi_n(a(n))}=1/\tilde{c}(n)$. It is clear that $\tilde{c}(n)$ and $c(n)$ are equivalent as $n$ tends to infinity. Therefore, $C_N$ has a threshold width of order $1/c(n)$. Furthermore, suppose now that $c$ is nondecreasing. Since $n\leq N(n)\leq 2n$, we have:
$$\forall n\in\NN^*,\;c(n)\leq c(N(n))\leq c(2n)\;.$$
Now, we suppose that $c(2n)=\Theta (c(n))$. Then $c(N(n))=\Theta (c(n))$. Finally, $C_N$ has a threshold width of order $1/c(N(n))$.
Whence the result.
{\hfill $\square$ \vskip 2mm \noindent} 

\section{How to get a sharp threshold from two coarse ones}
\label{coarsesharp}
When the localisation of the threshold of $B$ is not bounded away from 0 and
1, Proposition \ref{propprodtensseuil} is useless in describing the threshold width of $A\otimes
B$. Moreover, one can say that its conclusion is not valid any longer. Indeed,
consider the failure subset $B_k$ defined in Lemma \ref{lemseuillogn}:
$$B_k=A_{\lfloor\log_2 k\rfloor,\lfloor\log_2 k\rfloor}\otimes A_{1,\lfloor k/\log_2k\rfloor}\subset\{0,1\}^K\;.$$
One can easily compute the probability of $A_{1,n}$:
$$\mu_p(A_{1,n})=1-(1-p)^n\;,$$
therefore,
$$p(\eps)=1-(1-\eps)^{1/n}\;.$$
Thus, when  $n$ tends to infinity, $p(\eps)$ is equivalent to
$\log (1/(1-\eps))/n$, and $\tau(A_{1,n},\eps)$ is equivalent to $\log
((1-\eps)/\eps)/n$. This is the typical example of a coarse
threshold. Similarly, we get:
$$\mu_p(A_{n-1,n})=p^n\;,$$
therefore,
$$p(\eps)=\eps^{1/n}\;.$$
Thus, when  $n$ tends to infinity, $p(\eps)$ is equivalent to
$1-\log (1/\eps))/n$, and $\tau(A_{n-1,n},\eps)$ is again equivalent to $\log
((1-\eps)/\eps)/n$.\\
Consequently, $A_{\lfloor\log_2 k\rfloor,\lfloor\log_2 k\rfloor}$ and
$A_{1,\lfloor k/\log_2k\rfloor}$ have threshold widths of order, respectively
$1/ \log k$ and $\log k/k$. According to Lemma \ref{lemseuillogn}, their
product is a subset of $\{0,1\}^K$, where $K=\lfloor\log_2
k\rfloor\times\lfloor k/\log_2k\rfloor$, and has a threshold width of order
$1/\log K$. This is much bigger than the order $1/K$ which one would get if the
conclusion of Proposition \ref{propprodtensseuil} remained valid.\\
Nevertheless, this example witnesses an interesting phenomenon. The two subsets $A_{\lfloor\log_2 k\rfloor,\lfloor\log_2 k\rfloor}$ and
$A_{1,\lfloor k/\log_2k\rfloor}$ of $\{0,1\}^K$ clearly have coarse
thresholds, but their product has a sharp one. We shall prove this to be
a very general behaviour: as soon as $A$ and $B$ have thresholds
located respectively at 0 and 1, even if these are coarse, their product
$A\otimes B$ has a sharp threshold.\\
To this end, we shall use a well known tensorisation property of
the entropy. The major role of this property in concentration and threshold
topics has been pointed out many times (see for example
\cite{Ledoux96,Boucheronetal03,Rossignol05}). First, let us recall the definition of the entropy of a non negative function
$f$ on a probability space $(\mathcal{X},\mu)$:
$$\Ent_{\mu}(f)= \int f(x)\log f(x)\;d\mu(x)-\int f(x)\;d\mu(x)\log\int  f(x)\;d\mu(x)\;.$$
Entropy satisfies the following tensorisation inequality (see for instance
Ledoux \cite{Ledoux96}, Proposition 5.6 p.~98) : for every non
negative function $f$ on $\{0,1\}^n$,
\begin{equation}
\label{eqtens}
\sum_{i=1}^n\EE_{\mu_p}\left(\Ent_{x_i}(f)\right)\geq \Ent_{\mu_p}(f)\;,
\end{equation}
where $\Ent_{\mu_i}$ means that only the $i$-th coordinate is concerned with the
integration. The following lemma is the key towards the main result of this
section, Proposition \ref{propcoarsesharp}.
\begin{lem}
\label{lementropie}
Let $A$ be a monotone subset of $\{0,1\}^n$. Then, for every $p\in[0,1]$,
\begin{eqnarray}
\label{eqtensp} p\log
  \frac{1}{p}\frac{d\mu_p(A)}{dp}&\geq &\mu_p(A)\log\frac{1}{\mu_p(A)}\;,\\
\label{eqtens1-p} (1-p)\log
  \frac{1}{1-p}\frac{d\mu_p(A)}{dp}& \geq &(1-\mu_p(A))\log\frac{1}{1-\mu_p(A)}\;.
\end{eqnarray}
\end{lem}
\begin{proof}
The following formula is easy to obtain, by considering the derivative
of $\mu_p(x)$ with respect to $p$. For any real function $f$ on $\{0,1\}^n$,
\begin{equation}
\label{eqRussobis}
\frac{d}{dp}\int\!\!f(x) d\mu_p(x)=\sum_{i=1}^n\int\nabla_if(x)\;d\mu_p(x)\;,
\end{equation}
where
$$\forall x\in\{0,1\}^n,\;\nabla_if(x) =f(x_1,\ldots ,x_{i-1},1,x_{i+1},\ldots ,x_n)-f(x_1,\ldots ,x_{i-1},0,x_{i+1},\ldots ,x_n)\;.$$
On the other hand, if $A$ is a monotone subset,
\begin{eqnarray*}
\Ent_{x_i}(\II_A)&=&\int \II_A\log \frac{\II_A}{\int \II_A\;dx_i}\;dx_i\;,\\
&=&p\log\frac{1}{p}\nabla_i\II_A\;.
\end{eqnarray*}
Therefore, 
$$\sum_{i=1}^n\EE_{\mu_p}\left(\Ent_{x_i}(f)\right)=p\log \frac{1}{p}\frac{d\mu_p(A)}{dp}\;.$$
Remark that for any subset $A$,
$$\Ent_{\mu_p}(\II_A)=\mu_p(A)\log
\frac{1}{\mu_p(A)}\;.$$
Thus, when applied to $f=\II_{A}$, equation (\ref{eqtens}) gives
$$p\log \frac{1}{p}\frac{d\mu_p(A)}{dp}\geq \mu_p(A)\log
\frac{1}{\mu_p(A)}\;.$$
Now, remark that
\begin{eqnarray*}
\Ent_{x_i}(\II_{A^c})&=&\int \II_{A^c}\log \frac{\II_{A^c}}{\int \II_{A^c}\;dx_i}\;dx_i\;,\\
&=&p\log\frac{1}{p}\nabla_i\II_A\;.
\end{eqnarray*}
Since
$$\Ent_{\mu_p}(\II_{A^c})=\mu_p({A^c})\log
\frac{1}{\mu_p({A^c})}\;,$$
when applied to $f=\II_{A^c}$, equation \ref{eqtens} gives
$$(1-p)\log \frac{1}{1-p}\frac{d\mu_p(A)}{dp}\geq (1-\mu_p(A))\log
\frac{1}{1-\mu_p(A)}\;.$$
{\hfill $\square$ \noindent} \end{proof}

\begin{prop}
\label{propcoarsesharp}
Let $(r_n)_{n\in\NN}$ and $(m_n)_{n\in\NN}$
be two increasing sequences of integers, $A\subset\{0,1\}^{r_n}$
and $B\subset\{0,1\}^{m_n}$ be two non trivial monotone properties.
Suppose that the threshold of $A$ is located at 1, and the one of $B$ is
located at 0:
$$\forall \eps\in]0,1[,\;p_{A,\eps}\xrightarrow[n\rightarrow +\infty]{}1\;,$$
and
$$\forall \eps\in]0,1[,\;p_{B,\eps}\xrightarrow[n\rightarrow +\infty]{}0\;,$$
Then, the thresholds of $A\otimes B$ and $B\otimes A$ are sharp.
\end{prop}
\begin{proof}
First, let us remark that $A\otimes B$ has a coarse threshold if and only if,
for every $\eps\in]0,1[$ (see Friedgut and Kalai \cite{FriedgutKalai} or
Rossignol \cite{Rossthese}, p.99),
$$p_{A\otimes B,\eps}(1-p_{A\otimes B,\eps})\left.\frac{d\mu_p(A\otimes
    B)}{dp}\right|_{p=p_{A\otimes B,\eps}}\xrightarrow[n\rightarrow
    +\infty]{}+\infty\;.$$
Now, for all monotone subsets $A$ and $B$, and every $p$ in $]0,1[$, we get from
Proposition \ref{propprodtensgen}:
\begin{equation}
\label{eqprobatens}\mu_p(A\otimes B)  =\mu_{\mu_p(A)}(B)\;.
\end{equation}
Let us denote by $f(p)$ the quantity $\mu_p(A\otimes B)$. We have therefore:
$$f'(p)=\frac{d\mu_p(A)}{dp}\times\left.\frac{d\mu_q(B)}{dq}\right|_{q=\mu_p(A)}\;.$$
According to Lemma \ref{lementropie},
\begin{eqnarray*}
 \frac{d\mu_p(A)}{dp}&\geq &\frac{\mu_p(A)\log\frac{1}{\mu_p(A)}}{p\log
  \frac{1}{p}}\;,\\
 \left.\frac{d\mu_q(B)}{dq}\right|_{q=\mu_p(A)}& \geq &\frac{(1-f(p))\log\frac{1}{1-f(p)}}{(1-\mu_p(A))\log
  \frac{1}{1-\mu_p(A)}}\;.
\end{eqnarray*}
Thus,
\begin{equation}
\label{eqtensorparser}
f'(p)\geq \frac{(1-f(p))\log\frac{1}{1-f(p)}}{p\log
  \frac{1}{p}}\times
\frac{\mu_p(A)\log\frac{1}{\mu_p(A)}}{(1-\mu_p(A))\log\frac{1}{1-\mu_p(A)}}\;,
\end{equation}
We shall focus on $A\otimes B$, the threshold of $B\otimes A$ can be treated
in the same way, by switching the roles of $A$ and $B$. Now, suppose that the
threshold of $B$ is located at 0. Let $\eps$ belong to $]0,1[$. From
equation (\ref{eqprobatens}), one get
$$p_{A\otimes B,\eps}=p_{A,p_{B,\eps}}\;.$$
Let us note $p_\eps:=p_{A\otimes B,\eps}$. We get from inequality \ref{eqtensorparser}:
$$f'(p_\eps)\geq \frac{(1-\eps)\log\frac{1}{1-\eps}}{p_{\eps}\log
  \frac{1}{p_{\eps}}}\times
\frac{p_{B,\eps}\log\frac{1}{p_{B,\eps}}}{(1-p_{B,\eps})\log\frac{1}{1-p_{B,\eps}}}\;.
$$
Since $p_{B,\eps}$ tends to zero as $n$ tends to infinity,
\begin{eqnarray}
\nonumber f'(p_\eps)&\geq& \frac{(1-\eps)\log\frac{1}{1-\eps}}{p_{\eps}\log
  \frac{1}{p_{\eps}}}\times
\left(\log\frac{1}{p_{B,\eps}}+o(1)\right)\;,\\
\label{eqcoarsesharp1} p_\eps(1-p_\eps)f'(p_\eps)&\geq& (1-\eps)\log\frac{1}{1-\eps} \frac{(1-p_{\eps})}{\log
  \frac{1}{p_{\eps}}}\times
\left(\log\frac{1}{p_{B,\eps}}+o(1)\right)\;.
\end{eqnarray}
Since $\log 1/x$ is equivalent to $1-x$ as $x$ tends to 1, if $p_{\eps}$ is bounded away from zero, inequality (\ref{eqcoarsesharp1}) 
implies that $A\otimes B$ has a sharp threshold. If $p_{\eps}$ is not bounded
away from zero, one need to show that the location of $A$ at 1 implies that $\log 1/p_{\eps}$ is asymptotically negligible compared to
  $\log 1/p_{B,\eps}$. This follows from inequality (\ref{eqtensp}). Indeed,
  integrating this inequality between $p_{\eps}=p_{A,p_{B,\eps}}$ and
  $p_{A,\eps}$ leads to:
\begin{eqnarray*}
\left\lbrack\log\frac{1}{|\log
    \mu_p(A)|}\right\rbrack_{p_{A,p_{B,\eps}}}^{p_{A,\eps}}&\geq &\left\lbrack\log\frac{1}{|\log
    p|}\right\rbrack_{p_{A,p_{B,\eps}}}^{p_{A,\eps}}\;,\\
\log\frac{1}{|\log \eps|}-\log\frac{1}{|\log
    p_{B,\eps}|}&\geq & \log\frac{1}{|\log  p_{A,\eps}|}-\log\frac{1}{|\log
    p_{A,p_{B,\eps}}|}\;,\\
\log \frac{\log\frac{1}{p_{A,p_{B,\eps}}}}{\log\frac{1}{p_{B,\eps}}} & \leq & \log\frac{1}{\log
    1/\eps}+\log\frac{1}{|\log  p_{A,\eps}|}\;,\\
\frac{\log\frac{1}{p_{A,p_{B,\eps}}}}{\log\frac{1}{p_{B,\eps}}} & \leq &
    \frac{\log\frac{1}{ p_{A,\eps}}}{\log 1/\eps}\;.
\end{eqnarray*}
Since $p_{A,\eps}$ tends to 1 as $n$ tends to infinity, $\log 1/p_{A\otimes B,\eps}$ is asymptotically negligible compared to
  $\log 1/p_{B,\eps}$.
{\hfill $\square$ \noindent} \end{proof}


\bibliographystyle{apt}

\end{document}